\documentclass[a4paper,tbtags,draft]{amsart}
\usepackage{amssymb,latexsym}
\usepackage{amsxtra}
\usepackage{upref}
\usepackage[T1]{fontenc}
\usepackage{fix-cm}
\usepackage{fourier}
\usepackage[scaled]{beramono}
\usepackage{stmaryrd}
\usepackage{fix-stmaryrd}
\usepackage[all]{xy}
\usepackage{prooftree}
\usepackage{url}

\theoremstyle{definition}
\newtheorem{definition}{Definition}[section]

\theoremstyle{remark}

\theoremstyle{plain}

\DeclareMathOperator{\type}{:}

\DeclareMathOperator{\exI}{exI}

\DeclareMathOperator{\FV}{FV}
\DeclareMathOperator{\LTerms}{LTerms}

\DeclareMathOperator{\LambdaTypes}{\lambda Types}

\DeclareMathOperator{\Ax}{Ax}

\DeclareMathOperator{\Obj}{Obj}
\DeclareMathOperator{\Hom}{Hom}
\DeclareMathOperator{\id}{1}
\DeclareMathOperator{\inj}{\iota}
\DeclareMathOperator{\eval}{ev}

\DeclareMathOperator{\Transpose}{tr}

\newsavebox{\mysavebox}
\newdimen\myvskip
\newcommand{\directdisplay}[1][0pt]{%
  \par\ \par\sbox{\mysavebox}{\textit{Proof.}\ \ }%
  \myvskip\ht\mysavebox\advance\myvskip\dp\mysavebox%
  \advance\myvskip\baselineskip\advance\myvskip\parskip%
  \advance\myvskip 1.9\abovedisplayskip%
  \advance\myvskip #1\vskip -\myvskip}

\newdir{ >}{{}*!/-5pt/{>}}
\newdir{ <}{{}*!/-5pt/{<}}
\newdir{ V}{{}*!/-5pt/{\vee}}
\newdir{ A}{{}*!/-5pt/{\wedge}}

\begin{document}

\title{A Simplified Definition of Logically Distributive Categories}
\date{\today}

\author{Marco Benini}
\address{Dipartimento di Scienze Teoriche e Applicate\\
  Universit\`{a} degli Studi dell'Insubria\\
  via Maz\-zini 5, I-21100 Varese (VA), ITALY}
\email{marco.benini@uninsubria.it}
\urladdr{http://marcobenini.wordpress.com/}

\thanks{The author has been supported by the project \emph{Correctness
    by Construction} (CORCON), EU 7\textsuperscript{th} framework
  programme, grant n.~PIRSES-GA-2013-612638, and by the project
  \emph{Abstract Mathematics for Actual Computation: Hilbert's Program
    in the 21\textsuperscript{st} Century}, funded by the John
  Templeton's Foundation. Opinions expressed here and the
  responsibility for them are the author's only, not necessarily
  shared by the funding institutions.\\
  The author has also to thank the kind hospitality of
  prof.~H.~Ishihara and Dr.~T.~Nemoto when he was a visiting
  researcher at JAIST, where most of the results have been developed.}

\keywords{Category theory, logically distributive categories} 

\subjclass[2010]{Primary: 18A30; Secondary: 18A40, 03F55}

\begin{abstract}
  In~\cite{BE13} the notion of logically distributive category has
  been introduced to provide a sound and complete semantics to
  multi-sorted first-order logical theories based on intuitionistic
  logic. In this note, it will be shown that the definition of
  logically distributive category can be simplified by dropping two
  requirements, the one saying that product must distribute over sum,
  and the one saying that product must distribute over existential
  quantification. 
\end{abstract}

\maketitle

The purpose of this note is to show how one can drop two conditions in
the definition of logically distributive categories. We remind
from~\cite{BE13} the original definition:
\begin{definition}[Logically distributive
  category]\label{def:logicallydistributivecategory} 
  Fixed a $\lambda$-signature $\Sigma = \left\langle S, F, R,
    \Ax\right\rangle$, a category $\mathbb{C}$ together with a map
  $M\colon \LambdaTypes(\Sigma) \to \Obj \mathbb{C}$ is said to be
  \emph{logically distributive} if it satisfies the following seven
  conditions:
  \begin{enumerate}
  \item\label{cond:1} $\mathbb{C}$ has finite products;
  \item\label{cond:2} $\mathbb{C}$ has finite co-products;
  \item\label{cond:3} $\mathbb{C}$ has exponentiation;
  \item\label{cond:4} $\mathbb{C}$ is \emph{distributive}, i.e., for
    every $A, B, C \in \Obj \mathbb{C}$, the arrow $\Delta = [\id_A
    \times \inj_1,\id_A \times \inj_2]\colon (A \times B) + (A \times
    C) \to A \times (B + C)$ has an inverse, where $[\_,\_]$ is the
    co-universal arrow of the $(A \times B) + (A \times C)$
    co-product, $\_ \times \_$ is the product arrow, see~\cite{GOL06},
    $\id_A$ is the identity arrow on $A$, and $\inj_1\colon B \to B +
    C$, $\inj_2\colon C \to B + C$ are the canonical injections of the
    $B + C$ co-product.
  \end{enumerate}
  For every $s \in S$, $A \in \LambdaTypes(\Sigma)$, and $x \in V_s$,
  let $\Sigma_A(x\type s)\colon \LTerms(\Sigma)(s) \to \mathbb{C}$ be
  the functor from the discrete category $\LTerms(\Sigma)(s) = \left\{
    t\type s \mid t\type s \in \LTerms(\Sigma)\right\}$ defined by
  $t\type s \mapsto M(A[t/x])$.  Also, for every $s \in S$, $A \in
  \LambdaTypes(\Sigma)$, and $x \in V_s$, let $\mathbb{C}_{(\forall
    x\type s.\, A)}$ be the subcategory of $\mathbb{C}$ whose objects
  are the vertexes of the cones on $\Sigma_A(x\type s)$ such that they
  are of the form $M B$ for some $B \in \LambdaTypes(\Sigma)$ and
  $x\type s \not\in \FV(B)$. Moreover, the arrows of
  $\mathbb{C}_{(\forall x\type s.\, A)}$, apart identities, are the
  arrows in the category of cones over $\Sigma_A(x\type s)$ having the
  objects of $\mathbb{C}_{(\forall x\type s.\, A)}$ as domain and
  $M(\forall x\type s.\, A)$ as co-domain.  Finally, for every $s \in
  S$, $A \in \LambdaTypes(\Sigma)$, and $x \in V_s$, let
  $\mathbb{C}_{(\exists x\type s.\, A)}$ be the subcategory of
  $\mathbb{C}$ whose objects are the vertexes of the co-cones on
  $\Sigma_A(x\type s)$ such that they are of the form $M B$ for some
  $B \in \LambdaTypes(\Sigma)$ and $x\type s \not\in
  \FV(B)$. Moreover, the arrows of $\mathbb{C}_{(\exists x\type s.\,
    A)}$, apart identities, are the arrows in the category of co-cones
  over $\Sigma_A(x\type s)$ having the objects of
  $\mathbb{C}_{(\exists x\type s.\, A)}$ as co-domain and $M(\exists
  x\type s.\, A)$ as domain.

  \begin{enumerate}\setcounter{enumi}{4}
  \item\label{cond:5} All the subcategories $\mathbb{C}_{(\forall
      x\type s.\, A)}$ have terminal objects, and all the
    subcategories $\mathbb{C}_{(\exists x\type s.\, A)}$ have initial
    objects;
  \item\label{cond:6} The $M$ map is such that
    \begin{enumerate}
    \item $M(0) = 0$, the initial object of $\mathbb{C}$;
    \item $M(1) = 1$, the terminal object of $\mathbb{C}$;
    \item $M(A \times B) = M A \times M B$, the binary product in
      $\mathbb{C}$;
    \item $M(A + B) = M A + M B$, the binary co-product in
      $\mathbb{C}$;
    \item $M(A \to B) = M B^{M A}$, the exponential object in
      $\mathbb{C}$;
    \item $M(\forall x\type s.\, A)$ is the terminal object in the
      subcategory $\mathbb{C}_{(\forall x\type s.\, A)}$;
    \item $M(\exists x\type s.\, A)$ is the initial object in the
      subcategory $\mathbb{C}_{(\exists x\type s.\, A)}$;
    \end{enumerate}
  \item\label{cond:7} For every $x \in V_s$, $A, B \in
    \LambdaTypes(\Sigma)$ with $x\type s \not\in \FV(A)$, $M A \times
    M(\exists x\type s.\, B)$ is an object of $\mathbb{C}_{(\exists
      x\type s.\, A \times B)}$ since, if $\left( M(\exists x\type
      s.\, B), \left\{ \delta_t \right\}_{t\type s \in
        \LTerms(\Sigma)}\right)$ is a co-cone over $\Sigma_{B}(x\type
    s)$, and there is one by condition~(\ref{cond:5}), then $\left( M
      A \times M(\exists x\type s.\, B), \left\{ \id_{M A} \times
        \delta_t \right\}_{t\type s \in \LTerms(\Sigma)}\right)$ is a
    co-cone over $\Sigma_{A \times B}(x\type s)$. Thus, there is a
    unique arrow $!\colon M(\exists x\type s.\, A \times B) \to M A
    \times M(\exists x\type s.\, B)$ in $\mathbb{C}_{(\exists x\type
      s.\, A \times B)}$. Our last condition requires that the arrow
    $!$ has an inverse.
  \end{enumerate}
\end{definition}

Actually, since conditions~(\ref{cond:1}) and~(\ref{cond:3}) amount to
say that $\mathbb{C}$ is Cartesian closed\footnote{Some authors define
  Cartesian closed categories as those having all finite limits and
  exponentiation, and reserve the adjective ``weakly'' for the
  categories having just products and exponentiation. Since the result
  holds for both definition, we avoid the adjective.}, the product
functors $A \times \mathord{-} $ preserve co-limits, because it has
the exponential functor $(\mathord{-})^A$ as right adjoint --- a
standard check that can be found in any textbook, see, e.g.,
\cite{BOR94A}. So, by condition~(\ref{cond:2}), one can deduce the
existence of the arrow required in condition~(\ref{cond:4}). Actually,
it is more convenient to have a direct construction of the arrows
making the interpretations of $A \times (B + C)$ and $(A \times B) +
(A \times C)$ equivalent. This can be found in~\cite{BE14A}.

By this result, condition~(\ref{cond:4}) is redundant.

In~\cite{BE14A}, the natural bijection $\Hom(A \times B, C) \cong
\Hom(A, C^B)$ has been explicitly constructed. Specifically,
$\Transpose\colon \Hom(A \times B, C) \to \Hom(A, C^B)$, which maps an
arrow to its exponential transpose, is the isomorphism in one
direction, and $\theta$ works in the opposite direction, where
$\theta$ derives from a co-universal construction.

Turning to condition~(\ref{cond:7}), we want to show that it is
redundant, too. This amount to construct an inverse for the unique
arrow
\begin{align*}
  \alpha\colon M(\exists x\type s.\, A \times B) \to M A \times
  M(\exists x\type s.\, B)
\end{align*}
with $x\type s\not\in\FV A$.

The corresponding diagram in $\langle \mathbb{C}, M\rangle$ is
\begin{align}\label{diagram:1}
  \vcenter{\xymatrix{
    M(\exists x\type s.\, A \times B) \ar[rr]^{\alpha} & & M A \times
    M(\exists x\type s.\, B) \\
    & \left\{M A \times M\left(B[t/x]\right)\right\}_{t\type s}
    \ar[ul]^{\left\{\exI_t\right\}_{t\type s}}
    \ar[ur]_{\left\{q_t\right\}_{t\type s}} & } }
\end{align}
so that we have two co-cones and $\alpha$ makes the diagram to
commute, i.e., it is an arrow in the category of co-cones. Also, we
know that $M(\exists x\type s.\, A \times B)$ is the initial object in
the sub-category $\mathbb{C}_{\exists x\type s.\, A \times B}$, i.e.,
a co-limit.

We may think to the co-cone on the right of the
diagram~(\ref{diagram:1}) as the one produced by applying the product
functor $(A \times \mathord{-})$ to the following co-cone
\begin{align*}
  \xymatrix{
    M(\exists x\type s.\, B) \\
    \left\{M\left(B[t/s]\right)\right\}_{t\type s}
    \ar[u]_{\left\{\exI_t\right\}_{t\type s}} }
\end{align*}

Now, consider any co-cone of the form
\begin{align*}
  \xymatrix{
    M C \\
    \left\{M A \times M\left(B[t/s]\right)\right\}_{t\type s}
    \ar[u]_{\left\{p_t\right\}_{t\type s}} }
\end{align*}
where $x\type s \not\in\FV C$, i.e., any co-cone whose vertex lies in
$\mathbb{C}_{\exists x\type s.\, A \times B}$.

By exponentiation, we can construct
\begin{align*}
  \xymatrix{
    \left\{M A \times M\left(B[t/s]\right)\right\}_{t\type s}
    \ar[dr]^{\left\{p_t\right\}_{t\type s}}
    \ar@{.>}[dd]_{\left\{\id_{M A} \times \Transpose
        p_t\right\}_{t\type s}}& \\ 
    & M C \\
    M A \times {M C}^{M A} \ar[ur]_{\eval} }
\end{align*}

Since ${M C}^{M A} = M(A \supset C)$, $x\type s\not\in FV(A \supset
C)$, thus
\begin{align*}
  \xymatrix{
    {M C}^{M A} \\
    \left\{M\left(B[t/s]\right)\right\}_{t\type s}
    \ar[u]_{\left\{\Transpose p_t\right\}_{t\type s}} }
\end{align*}
is a co-cone whose vertex is in $\mathbb{C}_{\exists x\type s.\, B}$.

So, there is a co-universal arrow $\gamma$ as in the diagram
\begin{align*}
  \xymatrix{
    {M C}^{M A} & & M(\exists x\type s.\, B) \ar@{.>}[ll]_{\gamma} \\ 
    &\left\{M\left(B[t/s]\right)\right\}_{t\type s}
    \ar[ul]^{\left\{\Transpose p_t\right\}_{t\type s}}
    \ar[ur]_{\left\{\exI_t\right\}_{t\type s}} & } 
\end{align*}

Remembering that $\theta$ is the inverse of $\Transpose$, as
previously remarked, and applying the product functor, we get
\begin{align*}
  \xymatrix{
    M C & & M A \times M(\exists x\type s.\, B)
    \ar@{.>}[ll]_{\theta(\gamma)} \\ 
    & \left\{M A \times M\left(B[t/s]\right)\right\}_{t\type s}
    \ar[ul]^{\left\{p_t\right\}_{t\type s}}
    \ar[ur]_{\quad\left\{\id_{M A} \times \exI_t\right\}_{t\type s}} &
  }
\end{align*}

So, $M A \times M(\exists x\type s.\, B)$ is initial in
$\mathbb{C}_{\exists x\type s.\, A \times B}$ and thus isomorphic to
$M(\exists x\type s.\, A \times B)$ via the $\alpha$ arrow of
diagram~(\ref{diagram:1}) in one direction, and the $\theta(\gamma)$
arrow above in the opposite.

\bibliographystyle{amsalpha}
\bibliography{main}

\end{document}